\documentclass[12pt]{article}

\newtheorem{theorem}{Theorem}[section]
\newtheorem{lemma}[theorem]{Lemma}

\newtheorem{definition}[theorem]{Definition}

\newtheorem{remark}[theorem]{\sc Remark}

\newfam\msbfam
\font\tenmsb=msbm10  scaled \magstep1 \textfont\msbfam=\tenmsb
\font\sevenmsb=msbm7 scaled \magstep1 \scriptfont\msbfam=\sevenmsb
\font\fivemsb=msbm5  scaled \magstep1 \scriptscriptfont\msbfam=\fivemsb
\def\Bbb{\fam\msbfam \tenmsb}

\def\RR{{\Bbb R}}
\def\CC{{\Bbb C}}
\def\QQ{{\Bbb Q}}
\def\NN{{\Bbb N}}
\def\ZZ{{\Bbb Z}}
\def\II{{\Bbb I}}
\def\TT{{\Bbb T}}
\def\BB{{\Bbb B}}

\def\ss{\subseteq}
\def\ra{\rightarrow}

\def\O{\Omega}

 \def\HollowBox #1#2{{\dimen0=#1 \advance\dimen0 by -#2       
       \dimen1=#1 \advance\dimen1 by #2                       
        \vrule height #1 depth #2 width #2                    
        \vrule height 0pt depth #2 width #1                   
        \llap{\vrule height #1 depth -\dimen0 width \dimen1}%
       \hskip -#2                                             
       \vrule height #1 depth #2 width #2}}                   
 \def\BoxOpTwo{\mathord{\HollowBox{6pt}{.4pt}}\;}             

\def\endpf{\hfill $\BoxOpTwo$ \smallskip \\ }

\def\dbar{\overline{\partial}}

\newfam\msbbfam
\font\tenmsbb=msbm10  scaled \magstep1 \textfont\msbbfam=\tenmsbb
\font\sevenmsbb=msbm7  scaled \magstep1 \scriptfont\msbbfam=\sevenmsbb
\font\fivemsbb=msbm5    scaled \magstep1 \scriptscriptfont\msbbfam=\fivemsbb

\usepackage{graphicx}

\usepackage{amsmath}

\begin{document}

\begin{center}
\Large \bf The Greene-Krantz Conjecture in Dimension Two\footnote{{\bf Subject 
Classification Numbers:} 32M05, 32M12, 32M25  .}
\end{center}
\vspace*{.12in}

\begin{center}
Steven G. Krantz
\end{center}

\date{\today}

\begin{quote}
{\bf Abstract:}   We give a proof of the Greene-Krantz conjecture on convex domains
in $\CC^2$.  Curiously, the proof technique depends on subelliptic estimates
for the $\overline{\partial}$ problem.
\end{quote}
\vspace*{.25in}

\markboth{STEVEN G. KRANTZ}{THE BERGMAN PROJECTION}

\section{Introduction}

The last thirty-five years have seen a flourishing of the study
of the automorphism groups of smoothly bounded domains in $\CC^n$.
The subject has an unusual nature, because the only smoothly bounded
domain with transitive automorphism group is the unit ball $B$ (see [WON]).  So we
tend to instead focus our attention on the more general class of domains with {\it noncompact automorphism
group}.  It is a classical result of Cartan that such a domain $\Omega$ has
the property that there is a point $P \in \Omega$ and a point $X \in \partial \Omega$
and automorphisms (i.e., biholomorphic selfmaps of $\Omega$) $\varphi_j$ such
that $\varphi_j(P) \rightarrow X$ as $j \ra \infty$.   We call $X$ a {\it boundary
orbit accumulation point}.

Naturally we are interested in the geometric nature of the point $X$.  It is known (see [GRK1]) that $X$ must in fact
be a point of pseudoconvexity.  But we wish to know more about the Levi
geometry of $X$.   With this thought in mind, the following conjecture has
been formulated (see [GRK1]):
\begin{quote}
{\bf Greene-Krantz Conjecture:}   Let $\Omega$ be a smoothly bounded domain in $\CC^n$.
Suppose that $X \in \partial \Omega$ is a boundary orbit accumulation
point for the automorphism group action in the sense that there are
automorphisms $\varphi_j$ and a point $P \in \Omega$
such that $\varphi_j(P) \ra X$ as $j \ra \infty$.  Then $X$ is a point
of finite type in the sense of Kohn/D'Angelo/Catlin.
\end{quote}

This conjecture has been the object of intense study for the
past twenty years or more, and there are a number of
interesting partial results---see for instance [KIMS], [KIMK], [KIK1], [KIK2],
[KIK3]. In the present paper we prove this conjecture for
smoothly bounded convex domains in complex dimension two.

It is a pleasure to thank Harold Boas and Emil Straube for useful comments
and suggestions.

\section{Notation and Basic Ideas}

We take it that the reader is familiar with complex domains and with
pseudoconvexity.  See [KRA1] for background and details.  When the
ambient space has complex dimension two, there are two notions
of finite type, and they are as follows:

\begin{definition} \rm  A {\em first order commutator} of vector fields
is an expression of the form
$$
[L, M] \equiv L M - M L . 
$$
Note that the commutator is itself a vector field.

Inductively, an $m^{\rm th}$ order commutator is the commutator 
of an $(m - 1)^{\rm st}$ order commutator and a vector field $L$.
\end{definition}

\begin{definition} \rm  A {\it holomorphic vector field} is any linear combination
of the expressions 
$$
\frac{\partial}{\partial z_1} \quad , \quad  \frac{\partial}{\partial z_2} 
$$
with coefficients in the ring of $C^\infty$ functions.

A {\it conjugate holomorphic vector field} is any linear combination of the
expressions
$$
\frac{\partial}{\partial \overline{z}_1} \quad , \quad  \frac{\partial}{\partial \overline{z}_2} 
$$
with coefficients in the ring of $C^\infty$ functions.
\end{definition}

\begin{definition}  \rm
Let $M$ be a vector field defined on the boundary of $\O = \{z \in \CC^2: \rho(z) < 0\}.$  
We say
that $M$ is {\em tangential} if $M \rho = 0$ at each point of $\partial \Omega.$
\end{definition}

Now we define a gradation of vector fields which will be the basis for our
definition of analytic type.  Throughout this section 
$\O = \{z \in \CC^2: \rho(z) < 0\}$ and $\rho$ is $C^\infty$ with $\nabla \rho \ne 0$ on $\partial \Omega$.
If $X \in \partial \Omega$ then we may make a change
of coordinates so that $\partial \rho/\partial z_1 (X) \ne 0.$  Define
the holomorphic vector field
$$
   L = \frac{\partial \rho}{\partial z_2} \frac{\partial}{\partial z_1}  
        - \frac{\partial \rho}{\partial z_1} \frac{\partial}{\partial z_2} 
$$
and the conjugate holomorphic vector field
$$
  \overline{L} = \frac{\partial \rho}{\partial \overline{z}_2} \frac{\partial}{\partial \overline{z}_1} 
       - \frac{\partial \rho}{\partial \overline{z}_1} \frac{\partial}{\partial \overline{z}_2} . 
$$
Both $L$ and $\overline{L}$ are tangent to the boundary because $L \rho = 0$
and $\overline{L} \rho = 0.$  They are both non-vanishing near $X$ by our
normalization of coordinates.

The real and imaginary parts of $L$ (equivalently of $\overline{L}$) generate
(over the ground field $\RR$) 
the complex tangent space to $\partial \Omega$ at all points near $X$.
The vector field $L$ alone generates the space of all holomorphic
tangent vector fields and $\overline{L}$ alone generates the space of 
all conjugate holomorphic tangent vector fields.

\begin{definition} \rm  Let ${\cal L}_1$ denote the module, over the 
ring of $C^\infty$ functions, generated by $L$ and $\overline{L}.$
Inductively, ${\cal L}_\mu$ denotes the module generated
by ${\cal L}_{\mu - 1}$ and all commutators of the form
$[F,G]$ where $F \in {\cal L}_1$ and $G \in {\cal L}_{\mu - 1}.$
\end{definition}

Clearly ${\cal L}_1 \ss {\cal L}_2 \ss \cdots.$  Each ${\cal L}_\mu$ is closed under
conjugation.  {\em It is not generally the case that 
$\cup_\mu {\cal L}_\mu$ is the entire three-dimensional tangent space
at each point of the boundary.}  A counterexample is provided by
$$
\O = \{z \in \CC^2: |z_1|^2 + 2e^{-1/|z_2|^2} < 1\}  
$$
and the point $X = (1,0).$

\begin{definition} \rm  Let $\O = \{\rho < 0\}$ be a smoothly bounded domain in
$\CC^2$ and let $X \in \partial \Omega.$  
We say that $\partial \Omega$ is of {\em finite analytic type $m$ at $X$}
if $\langle \partial \rho(X), F(X) \rangle = 0$ for
all $F \in {\cal L}_{m-1}$ while $\langle \partial \rho(X), G(X) \rangle \ne 0$
for some $G \in {\cal L}_{m}.$
In this circumstance we call $X$ a {\it point of analytic type $m$}.
\end{definition}
\vskip.1in

Now we turn to a precise definition of finite geometric type.
Let $D$ denote the unit disc in the complex plane.
If $X$ is a point in the boundary of a smoothly bounded domain
then we say that an analytic disc $\phi: D \ra \CC^2$ is 
a {\em non-singular disc tangent to $\partial \Omega$ at $X$} if
$\phi(0) = X, \phi'(0) \ne 0,$ and $(\rho\circ \phi)'(0) = 0.$

\begin{definition} \rm
Let $\O = \{\rho < 0\}$ be a smoothly bounded domain and
$X \in \partial \Omega.$  Let $m$ be a non-negative integer.
We say that $\partial \Omega$ is of finite geometric
type $m$ at $X$ if the following condition holds:  
there is a non-singular disc $\phi$ 
tangent to $\partial \Omega$ at $X$
such that, for small $\zeta,$ 
$$
|\rho \circ \phi(\zeta)| \leq C |\zeta|^m 
$$
BUT there is no non-singular disc $\psi$
tangent to $\partial \Omega$ at $X$ such that, for small $\zeta,$
$$
|\rho \circ \phi(\zeta)| \leq C |\zeta|^{(m+1)} . 
$$
In this circumstance we call $X$ a point of finite geometric
type $m.$
\end{definition}

The principal result about finite type in dimension two is the
following theorem (see [KRA1, \S 11.5]):

\begin{theorem}  \sl 
Let $\O = \{\rho < 0\} \ss \CC^2$ be smoothly bounded and
$X \in \partial \Omega.$  The point $X$ is of finite geometric type
$m \geq 2$ if and only if it is of finite analytic type $m.$
\end{theorem}

Now let us say a few words about subelliptic estimates.   A partial
differential operator ${\cal L}$ of order $k$ is said to satisfy {\it elliptic
estimates} if, whenever ${\cal L} u = f$ and $f$ lies in
the Sobolev space $W^s$ then $u$ lies in the Sobolev space $W^{s + k}$.
The operator is said to satisfy {\it subelliptic estimates} if the index
$s + k$ in the conclusion is replaced by $s + k'$ for some $0 < k' < k$.
The $\overline{\partial}$-Neumann operator on a strongly pseudoconvex domain,
and more generally on a finite type domain, is known to satisfy a subelliptic
(but definitely not an elliptic) estimate.  See [CAT1]--[CAT2], [KRA3], [FOK] for the
details.   It is also possible to express the subellipticity condition in
terms of Lipschitz or Besov spaces rather than Sobolev spaces.  We leave the details
for the interested reader.

\section{The Main Argument}

The result that we shall actually prove in this paper is the following:

\begin{theorem} \sl 
Let $\Omega$ be a smoothly bounded, convex
domain in $\CC^n$. Suppose that $X \in \partial \Omega$ is a
boundary orbit accumulation point for the automorphism group
action in the sense that there are automorphisms $\varphi_j$, a point $P \in \Omega$, and
a point $X \in \Omega$ such that 
$\varphi_j(P) \ra X$ as $j \ra \infty$.  Then $X$ is a point of finite type
in the sense of Kohn/D'Angelo/Catlin. 
\end{theorem}

Now fix a smoothly bounded domain $\Omega \ss \CC^2$.  Assume that
$P \in \Omega$ and $X \in \partial \Omega$ and that there are
automorphisms $\varphi_j$ of $\Omega$ such that $\varphi_j(P) \ra X$
as $j \ra \infty$.  Note that, because the domain $\Omega$ is smoothly bounded
and pseudoconvex, it is complete in the Bergman metric (see [OHS]).

Now consider a small Bergman metric ball $\beta$ centered at $P$.   Choose $j_1$ so that
$\beta_1 \equiv \varphi_{j_1}(\beta)$ is disjoint from $\beta$, and so that the Euclidean distance
of $\beta_1$ to the boundary is about $2^{-1}$.  Now choose $j_2$ 
so that $\beta_2 \equiv \varphi_{j_2}(\beta)$ is disjoint from $\beta$ and $\varphi_1(\beta)$ and
so that the Euclidean distance of $\beta_2$ to the boundary is about $2^{-2}$.  Keep going.

Now fix a $\overline{\partial}$-closed $(0,1)$ form $\psi$ with
$C_c^\infty$ coefficients that is supported in $\beta$. Define $\psi_\ell =
(\varphi_{j_\ell}^{-1})^* \psi$. Thus $\psi_\ell$ is a $\overline{\partial}$-closed
$(0,1)$ form with $C_c^\infty$ coefficients supported on $\beta_\ell$.  Because
of the derivative of $\varphi_{j_\ell}$, $\psi_\ell$ has supremum norm about $2^{-\ell}$.
That will mean that the sum of the $\psi_\ell$ will
have an $L^2$ or Sobolev norm that converges.

If we write $\psi_\ell = \psi_\ell^1 d\overline{z}_1 + \psi_\ell^2 \overline{z}_2$, then 
we may note that the equation $\overline{\partial} u_\ell = \psi_\ell$ can
be solved with one of the simple equations 
$$
u_\ell^1(z_1, z_2) = \mathop{\int \!\!\!\int}_{\zeta \in \CC} \frac{\psi_\ell^1(\zeta, z_2)}{\zeta - z_1} \, dA(\zeta)
$$
or
$$
u_\ell^2(z_1, z_2) = \mathop{\int \!\!\!\int}_{\zeta \in \CC} \frac{\psi_\ell^2(z_1, \zeta)}{\zeta - z_2} \, dA(\zeta)	  \, ,
$$
see [KRA1, \S 1.1].   And it turns out that $u_\ell^1 = u_\ell^2$.  

It follows from standard results on fractional integration (see [STE]) that,
if $\psi_\ell$ is in some Sobolev class $W^s$ then $u_\ell^m$ will be in a smoother
Sobolev class $W^{s'}$, with $s' > s$, in the $m^{\rm th}$ variable, $m = 1, 2$.  And now
a simple argument with the triangle inequality shows that $u_\ell^m$ lies
in $W^{s''}$ as a function of both variables for some $s' \geq s'' > s$, $m = 1, 2$.  So we see that the $\overline{\partial}$ problem
satisfies a subelliptic estimate on $\psi_\ell$.  

But in fact, thanks to the intervention of the automorphisms $\varphi_\ell$, the $\overline{\partial}$
problem satisfies the {\it very same} subelliptic estimate for each $\psi_\ell$.
As a result, the $\overline{\partial}$ problem satisfies a subelliptic estimate on
the form
$$
\psi \equiv \sum_\ell \psi_\ell \, .
$$

Now it is definitely not the case that the $\overline{\partial}$-closed (0,1) forms with
$C_c^\infty$ coefficients are dense in any space of forms with Sobolev coefficients.
But we shall be able to argue that they {\it are} dense in certain forms that we
care about.  See also the footnote below.

We have the following lemma:

\begin{lemma} \sl
If the boundary orbit accumulation point $X$ is of infinite type, then for each $\epsilon > 0$ there is
a $\overline{\partial}$-closed $(0,1)$ form $f$ on $\Omega$ with $L^2$ coefficients so
that the equation $\overline{\partial} u = f$ does {\it not} have any solution
in the Besov space of order $\epsilon > 0$.
\end{lemma}
{\bf Proof:}  The idea for the proof goes back to an old result of
Kerzman (see [KER]) and is reasonably well known.  See also [KRA1, \S 10.3].
We sketch the idea here.

We may assume that $X = (1,0) \in \partial \Omega$ and that the complex normal direction at $X$ is
$\langle 1, 0\rangle$.  With these normalizations, we define
$$
f = \frac{d\overline{z}_2}{\log(1 - z_1)} \, .
$$
By the convexity of $\Omega$, it is clear that the principal branch of the logarithm 
is well defined and that $f$ has bounded coefficients.

Now any solution of the equation $\overline{\partial} u = f$ will have the 
form
$$
u(z) = \frac{\overline{z}_2}{\log(1 - z_1)} + h(z_1, z_2) \, ,
$$ 
where $h$ is some holomorphic function on $\Omega$.  

Since $X$ is a point of infinite type then we know that, for any positive integer $m$, there is a nonsingular
complex curve $\mu_m: D \ra \CC^2$ that is tangent to order $2m$ with $\partial \Omega$ at $X$.  Let $\nu_X$
denote the Euclidean outward unit normal vector to $\partial \Omega$ at $X$.  Then, for $\delta > 0$ small, 
the analytic disc
$$
\biggl \{\mu_m(\zeta) - \delta \nu_X: |\zeta| < C \delta^{1/(2m)}, \zeta \in D \biggr \}
$$
lies in $\Omega$ (see [KRA2] for the elementary calculations needed to justify
this assertion).  Thus
$$
\theta_\delta: t \longmapsto \mu_m (C \delta^{1/(2m)} e^{it}) - \delta \nu_X \ \ , \quad 0 \leq t < 2\pi \ ,
$$
describes the boundary of an analytic disc in $\Omega$.

With $(1,0) \in \partial \Omega$ as our point of infinite type, let $m$ be a positive integer as above
and (by the well-known semicontinuity of type---see [KRA1, \S 11.5]) choose a neighborhood $W$ of
$(1,0)$ so that boundary points in $W$ are of finite type at least $2m$.  We may assume that $W$ is a
Euclidean ball, and that it lies in a tubular neighborhood of $\partial \Omega$.  Pick $\delta > 0$ small 
(small enough so the $\delta^{1/(2m)}$ is much less than the radius of $W$) and define
$$
\widetilde{\Omega} = (W \cap \Omega) \bigcup \biggl \{ z \in \Omega: \hbox{dist} (z, \partial \Omega) > \delta^{1/(2m)} \biggr \} \, .
$$

We examine
the complex line integral
$$
F(\delta, \zeta) = \oint_{\theta_\delta} u(\zeta_1 - 2\delta, \zeta_2 + z_2) - u(\zeta_1 - \delta, \zeta_2 + z_2) dz_2 
$$
for $\zeta \in \widetilde{\Omega}$.
We note that the curves
$$
t \longmapsto - \delta\nu_X + \mu(C' \delta^{1/(2m)} e^{it}) \qquad \hbox{and} \qquad t \longmapsto - 2\delta\nu_X + \mu(C' \delta^{1/(2m)} e^{it})
$$
both lie in $\Omega$ precisely because $X$ is a point of infinite type (more precisely, a point of type at least $2m$).  

Seeking a contradiction, if $u$ satisfies a Besov condition of order $\epsilon$, then we may straightforwardly 
estimate that
\begin{eqnarray*}
\|F(\delta)\|_{L^2(\zeta)} & \leq & \int_{\widetilde{\Omega}} \left | \oint_{\theta_\delta} u(\zeta_1 - 2\delta, \zeta_2 + z_2) - u(\zeta_1 - \delta, \zeta_2 + z_2) dz_2 \right |^2 dV(\zeta)^{/12}  \\
                           & \leq & \int_{\theta_\delta} \int_{\widetilde{\Omega}} |u(\zeta_1 - 2\delta, \zeta_2 + z_2) - u( \zeta_1 - \delta, \zeta_2 + z_2) |^2 dV(\zeta)^{1/2} d|z_2| \\
			   & \leq & \int_{\theta_\delta} \delta^\epsilon \, d|z_2| \\
			   & \approx & \delta^{\epsilon + 1/(2m)} \, .
\end{eqnarray*}

On the other hand,
\begin{eqnarray*}
  F(\delta, \zeta) & = & \int_{\theta_\delta} \frac{\overline{z}_2 + \overline{\zeta}_2}{\log (1 - \zeta_1 + 2\delta)} - 
			      \frac{\overline{z}_2 + \overline{\zeta}_2}{\log (1 - \zeta_1 + \delta)} \, dz	 \\
	    & = & \frac{\delta^{2/(2m)}}{\log(1 - \zeta_1 + 2\delta)} - \frac{\delta^{2/(2m)}}{\log(1 - \zeta_1 + \delta)}   \\
	    & \approx & C \cdot \frac{\delta^{1/m}}{\log^2( - \delta)} \, .	 \\
\end{eqnarray*}
As a result,
$$
\int_{\widetilde{\Omega}} |F(\delta)|^2 \, dV(\zeta)^{1/2} \approx \frac{\delta^{1/m}}{\log^2(- \delta)} \, .
$$

Comparing our two estimates, we find that
$$
\frac{\delta^{1/m}}{\log^2(-\delta)} \leq C \cdot \delta^{\epsilon + 1/(2m)}
$$
or
$$
\frac{\delta^{1/(2m)}}{\log^2(-\delta)} \leq C \cdot \delta^\epsilon \, .
$$
This is false as soon as $m \in \NN$ is large enough.
\endpf
\smallskip \\

The lemma tells us that, in the Besov space topology, the $\overline{\partial}$ problem does {\it not}
satisfy a subelliptic estimate.  But it is not difficult to see that the form
$$
f(z) = \frac{d\overline{z}_2}{\log(1 - z_1)}
$$
is the limit of forms with compact support.\footnote{And notice that, if $\psi_0$ is a $\overline{\partial}$-closed
$(0,1)$ form with $C_c^\infty$ coefficients on $\beta$, then we can consider the form 
$\psi$ on the union of $\beta$, $\varphi_{j_1}(\beta)$, etc. as described above and we can also
consider the ``shifted'' form $\tau$ given by $(\varphi_{j_1}^{-1})^*\psi$ on $\varphi_{j_1}(\beta)$, $(\varphi_{j_2}^{-1})^*\psi$ on $\varphi_{j_2}(\beta)$ (with
intervening automorphism $\varphi_{j_2} \circ \varphi_{j_1}^{-1}$), and so forth.  Then the difference
of these two forms is a $C_c^\infty$ form supported on $\beta$ alone.   So our arguments and estimates also apply to
forms that have compact support and are smooth.}  
For let $\rho_2$ be a $C_c^\infty$ function
that approximates $1/\log(1 - z_1)$ in the $L^2$ topology.   Now the formula
$$
v(z_1, z_2) = \mathop{\int\!\!\!\int} \frac{\rho_2(z_1, \zeta)}{\zeta - z_2} \, d\zeta
$$
satisfies 
$$
\frac{\partial}{\partial \overline{z}_2} v = \rho_2 \, .
$$
Note that (see [KRA1, \S 1.1]) $v \in C_c^\infty(\Omega)$.
Hence
$$
\rho_1(z) \equiv \frac{\partial}{\partial \overline{z}_1} v
$$
will give a form 
$$
R = \rho_1 d\overline{z}_1 + \rho_2 d\overline{z}_2
$$
that is $\overline{\partial}$-closed with $C_c^\infty$ coefficients.  And of 
course $R$ will approximate $f$ in the $L^2$ topology.

This approximation implies of course that the problem $\overline{\partial}u = f$, with $f$
as in the lemma, satisfies a subelliptic estimate in the Sobolev topology.  But that implies
that it satisfies a subelliptic estimate in the Besov topology.  And we have established
in the lemma that that is impossible.

We have proved that the boundary orbit accumulation point $X$ {\it cannot} be of infinite type.

\begin{remark} \rm
It is worth noting that the construction presented here---of the ball $\beta$ and subsequent
target balls $\varphi_{j_1}(\beta)$, $\varphi_{j_2}(\beta)$, etc., does not work
when the automorphism group is compact.   For, when the automorphism group
is compact, then these balls will no longer be pairwise disjoint.  Also the norms of
the $(\varphi_{j_k}^{-1})^* \psi$ will no longer vanish rapidly, so that the
series which is obtained by adding the forms supported on the different
balls will no longer converge.
\end{remark}

\section{Concluding Remarks}

In this paper we certainly have not proved the full Greene-Krantz
conjecture.  But we have proved a notable and interesting special
case.

There is certainly interest in developing techniques for attacking
the full conjecture, and we intend to attack that problem in future
papers.

\newpage

\noindent {\Large \sc References}
\vspace*{.2in}

\begin{enumerate}

\newfam\msbfam
\font\tenmsb=msbm10  scaled \magstep1 \textfont\msbfam=\tenmsb
\font\sevenmsb=msbm7 scaled \magstep1 \scriptfont\msbfam=\sevenmsb
\font\fivemsb=msbm5  scaled \magstep1 \scriptscriptfont\msbfam=\fivemsb
\def\Bbb{\fam\msbfam \tenmsb}

\def\RR{{\Bbb R}}
\def\CC{{\Bbb C}}
\def\QQ{{\Bbb Q}}
\def\NN{{\Bbb N}}
\def\ZZ{{\Bbb Z}}
\def\II{{\Bbb I}}
\def\TT{{\Bbb T}}
\def\BB{{\Bbb B}}

\item[{\bf [CAT1]}] D. Catlin, Necessary conditions for
subellipticity of the $\overline{\partial}-$Neumann problem,
{\it Ann. Math.} 117(1983), 147-172.

\item[{\bf [CAT2]}] D. Catlin, Subelliptic estimates for the
$\dbar$Neumann problem, {\it Ann. Math.} 126(1987), 131-192.

\item[{\bf [FOK]}]  G. B. Folland and J. J. Kohn, {\it The Neumann
Problem for the Cauchy-Riemann Complex}, Princeton University
Press, Princeton, NJ, 1972.

\item[{\bf [GRK1]}]  R. E. Greene and S. G. Krantz, Invariants
of Bergman geometry and results concerning the automorphism
groups of domains in $\CC^n$, {\it Geometrical and Algebraical
Aspects in Several Complex Variables} (Cetraro, 1989),
107--136, Sem.\ Conf., 8, EditEl, Rende, 1991.

\item[{\bf [KIMS]}]  S.-Y. Kim, Domains with hyperbolic orbit
accumulation boundary points, {\it J. Geom.\ Anal.} 22 (2012),
90--106.

\item[{\bf [KIMK]}] K.-T. Domains in $\CC^n$ with a piecewise
Levi flat boundary which possess a noncompact automorphism
group, {\it Math.\ Ann.} 292(1992), 575--586. ,

\item[{\bf [KIK1]}]  K.-T. Kim and S. G. Krantz, Complex scaling
and domains with non-compact automorphism group, {\it Illinois
Journal of Math.} 45(2001), 1273--1299.
			
\item[{\bf [KIK2]}]  K.-T. Kim and S. G. Krantz, Some new
results on domains in complex space with non-compact
automorphism group, {\it J.\ Math.\ Anal.\ Appl.} 281(2003),

\item[{\bf [KIK3]}] K.-T. Kim and S. G. Krantz, Complex scaling
and geometric analysis of several variables, {\it Bull.\
Korean Math.\ Soc.} 45(2008), 523--561.
		 
\item[{\bf [KRA1]}]  S. G. Krantz, {\it Function Theory of
Several Complex Variables}, 2nd. ed., American Mathematical
Society, Providence, RI, 2001.

\item[{\bf [KRA3]}]  S. G. Krantz, {\it Partial Differential
Equations and Complex Analysis}, CRC Press, Boca Raton, FL, 1992.

\item[{\bf [KRA2]}] S. G. Krantz, Characterizations of various
domains of holomorphy via $\dbar$-estimates and applications
to a problem of Kohn, {\it Illinois J.\ Math.} 23(1979),
267--285.

\item[{\bf [OHS]}] T. Ohsawa, A remark on the completeness of
the Bergman metric, {\it Proc.\ Japan Acad.} Ser.\ A Math.
Sci. 57(1981), 238--240.

\item[{\bf [WON]}] B. Wong, Characterizations of the ball in
$\CC^n$ by its automorphism group, {\em Invent. Math.}
41(1977), 253-257.

\end{enumerate}
\vspace*{.25in}

\begin{quote}
Steven G. Krantz  \\
Department of Mathematics \\
Washington University in St. Louis \\
St.\ Louis, Missouri 63130  \\
{\tt sk@math.wustl.edu}
\end{quote}

\end{document}